\documentclass[12pt]{article}
\usepackage{amsmath,amssymb,amsfonts}
\usepackage[utf8]{inputenc}
\usepackage{color}

\begin{document}

\begin{center}
\Large\bf On the convergence of formal Dulac series \\ satisfying an algebraic ODE
\end{center}

\begin{center}

Irina Goryuchkina\footnote{Keldysh Institute of Applied Mathematics of RAS, Russia, Moscow,\; igoryuchkina@gmail.com},\;
Renat Gontsov\footnote{Institute for Information Transmission Problems of RAS, Russia, Moscow, \; gontsovrr@gmail.com}	

\end{center}

\bigskip

\begin{abstract}

We propose a sufficient condition of the convergence of a Dulac series formally satisfying an algebraic ordinary
differential equation (ODE). Such formal solutions of algebraic ODEs appear rather often, in particular, the 
third, fifth, and sixth Painlev\'e equations possess formal Dulac series solutions, whose convergence follows 
from the proposed sufficient condition. 

\end{abstract}

\section{Introduction}

We consider an $n$-th order ODE
\begin{equation}\label{ODE}
F(x,y,\delta y,{\dots},\delta^n y)=0,
\end{equation}
where $F=F(x,y_0,y_1,\ldots,y_n)$ is a polynomial of $n+2$ complex variables and $\delta$ is the derivation
$x(d/dx)$. Assume (\ref{ODE}) has a formal {\it Dulac series} solution $\varphi$ of the form
$$
\varphi=\sum_{k=0}^{\infty}p_k(\ln x)\,x^{k}, \qquad p_k\in{\mathbb C}[t].
$$

Such series appeared in the 1920's in the works of Henry Dulac \cite{Du} on limit cycles of a planar vector
field as asymptotic expansions for the monodromy map (first return map) in a
neighbourhood of a polycycle. (More precisely, these were the series with more general powers $x^{\lambda_k}$ rather than $x^k$, where $\lambda_k$ formed a sequence of real numbers increasing to infinity.) Further, in the 1980's, Dulac series played an important role in finishing proofs of the finiteness theorems for limit cycles (see \cite{Il1} or \cite{Il2}).

Dulac series also appear as formal solutions of algebraic ODEs (of the Abel equations, Emden-Fowler type equations, Painlev\'e
equations, {\it etc.}) and this is a subject of our interest in the present work. First question that one may naturally pose in this context, is that of the convergence of such formal solutions. We propose the following sufficient condition of convergence.
\medskip

{\bf Theorem 1.} {\it Let the series $\varphi$ formally satisfy the equation under consideration:
$$F(x,\Phi):=F(x,\varphi,\delta\varphi,\dots,\delta^n \varphi)=0,$$
and let for each $j=0,\dots,n$ one have
$$
\frac{\partial F}{\partial y_j}(x,\Phi)=a_j x^m+b_j(\ln x) x^{m+1}+\dots,
$$
where $a_j\in\mathbb{C}$, $m\in{\mathbb Z}_+$ is the same for all $j$, and $b_j\in\mathbb{C}[t]$.
If $a_n\neq 0$ then for any open sector $S$ of sufficiently small radius, with the vertex at the origin and of the opening less than $2\pi$, the series $\varphi$ converges uniformly in $S$.}
\medskip

For example, the third, fifth and sixth Painlev\'e equations have formal solutions in Dulac series. S.\,Shimomura \cite{Sh1}, \cite{Sh2} has proved their convergence for the Painlev\'e V and VI using the connection of these equations with isomonodromic deformations of linear differential systems. One can also apply Theorem 1 to prove the convergence of formal Dulac series solutions of {\it all} the Painlev\'e equations (see examples at the last section of the paper).

Note that in the case of all $p_k=\rm const$, one has a formal {\it power series} solution $\varphi$ of (\ref{ODE}) and Theorem 1 becomes a well known sufficient condition of the convergence of such a solution obtained by B.\,Malgrange \cite{Ma}. We also note that combining the technique of the present paper with that of \cite{GG} on the convergence of {\it generalized} power series solutions of (\ref{ODE}), one can obtain a theorem similar to Theorem 1 for formal Dulac series of a more general form (with $x^{\lambda_k}$, $\lambda_k\in\mathbb{C}$, instead of $x^k$).

The paper is organized as follows. Theorem 1 is finally proved in Section 5, which is preceeded by some auxiliary statements: in Section 2 we pass from the initial ODE to a reduced one, in Section 3 there is proposed some auxiliary linear algebra, and in Section 4 we construct an ODE which is majorant for the reduced ODE. In the last Section 6 we give two examples of the
application of Theorem 1.

\section{An equation in the reduced form}

{\bf Lemma 1.} {\it Under the conditions of Theorem 1, there is $\ell'\in{\mathbb Z}_+$ such that for any $\ell\geqslant\ell'$
the transformation $$y=\sum\limits_{k=0}^{\ell}p_k(\ln x)x^k+x^\ell\, u$$ reduces the initial equation $(\ref{ODE})$ to the equation 
\begin{equation}\label{rODE}
L(\delta)u=x\,M(x, \ln x, u,\delta u,\dots,\delta^n u),
\end{equation}
where $$L(\delta)=\sum\limits_{j=0}^n a_j(\delta+\ell)^j, \quad a_n\ne0,$$ and 
$M$ is a polynomial of $n+3$ variables. Moreover, the polynomial $L$ does not vanish in the open right-half plane.}
\medskip

{\bf Proof.} The method of the Lemma 1 proof is standard and similar to that of the proof of the corresponding reduction lemma for power series
solutions (see \cite{Ma}).

For each non-negative integer $\ell$ the formal Dulac series $\varphi$ can be represented in the form
$$
\varphi=\sum\limits_{k=0}^{\ell}p_k(\ln x)x^k+x^{\ell}\sum\limits_{k=1}^{\infty}p_{k+\ell}(\ln x)x^k=:\varphi_{\ell}+x^{\ell}\psi,
$$
then
$$
\Phi=(\varphi,\delta\varphi,\ldots,\delta^n\varphi)=\Phi_{\ell}+x^{\ell}\Psi,
$$
where $\Phi_{\ell}=(\varphi_{\ell},\delta\varphi_{\ell},\ldots,\delta^n\varphi_{\ell})$ and $\Psi=(\psi,(\delta+\ell)\psi,\ldots,(\delta+\ell)^n\psi)$.
Applying Taylor's formula one obtains
\begin{eqnarray}\label{Taylor}
0=F(x,\Phi_{\ell}+x^{\ell}\Psi)=F(x,\Phi_{\ell})+x^{\ell}\sum_{j=0}^n\frac{\partial F}{\partial y_j}(x,\Phi_{\ell})\psi_j+ \nonumber \\
+\frac{x^{2\ell}}2\sum_{i,j=0}^n\frac{\partial^2 F}{\partial y_i\partial y_j}(x,\Phi_{\ell})\psi_i\psi_j+\ldots,
\end{eqnarray}
where $\psi_j=(\delta+\ell)^j\psi$.

Let us choose the number $\ell$ in such a way that the following two conditions hold:
$$
1)\; \ell>m, \qquad 2)\; L(\xi)=\sum_{j=0}^n a_j(\xi+\ell)^j\ne0 \quad\forall\xi\in\{{\rm Re}\,\xi>0\}
$$
(recall that the non-negative integer $m$ comes from the condition of Theorem 1).
\medskip

{\bf Definition.} Define the valuation of an arbitrary Dulac series $\varphi=\sum_{k=0}^{\infty}p_k(\ln x)x^k$
as
$$
{\rm val}(\varphi):=\min\{k\mid p_k\not\equiv0\}.
$$

Since Taylor's formula gives
$$
\frac{\partial F}{\partial y_j}(x,\Phi)-\frac{\partial F}{\partial y_j}(x,\Phi_{\ell})=x^{\ell}\sum_{i=0}^n
\frac{\partial^2 F}{\partial y_i\partial y_j}(x,\Phi_{\ell})\psi_i+\ldots,
$$
and ${\rm val}(\psi_i)\geqslant1$ for any $i$, one has
$$
\frac{\partial F}{\partial y_j}(x,\Phi_{\ell})=a_j x^m+\tilde b_j(\ln x) x^{m+1}+\dots, \quad \tilde b_j\in\mathbb{C}[t],
$$
for each $j=0,1,\ldots,n$, that is, the leading coefficient $a_j$ is preserved when one substitutes in $\frac{\partial F}{\partial y_j}$
the finite sum $\Phi_{\ell}$ instead of $\Phi$. Now, the relation (\ref{Taylor}) implies that
$$
{\rm val}(F(x,\Phi_{\ell}))\geqslant m+\ell+1.
$$
Dividing the relation (\ref{Taylor}) by $x^{m+\ell}$ one obtains the equation of the prescribed form (\ref{rODE}) with the formal Dulac series solution
$$
\psi=\sum\limits_{k=1}^{\infty}p_{k+\ell}(\ln x)x^k=\sum\limits_{k=1}^{\infty}P_k(\ln x)x^k.
$$
The lemma is proved.{\hfill $\Box$}
\medskip

{\bf Lemma 2.} {\it The formal Dulac series solution $\psi$ of $(\ref{rODE})$ is uniquely determined and the degree $\nu_k$ of each
polynomial $P_k$ satisfies the estimate $\nu_k\leqslant k\,C$, where $C$ is the degree of the polynomial $M$ with respect to $\ln x$.}
\medskip

{\bf Proof.} First, note that one has the following differentiation rule:
$$
\delta: \; P_k(\ln x)\,x^k \;\mapsto\; x^k\,\Bigl(k+\frac{d}{dt}\Bigr)P_k(t)|_{t=\ln x},
$$
hence
\begin{eqnarray*}
(\delta+\ell)^j: & & P_k(\ln x)\,x^k \;\mapsto\; x^k\,\Bigl(k+\ell+\frac{d}{dt}\Bigr)^jP_k(t)|_{t=\ln x},\\
L(\delta):  & & P_k(\ln x)\,x^k \;\mapsto\; x^k\,L\Bigl(k+\frac{d}{dt}\Bigr)P_k(t)|_{t=\ln x}.
\end{eqnarray*}
%
Thus, substituting $\psi=\sum_{k=1}^{\infty}P_k(\ln x)x^k$ into the equation
$$
L(\delta)u=x\,M(x,\ln x, u,\delta u,\dots,\delta^n u)
$$
we obtain the relation whose both sides are Dulac series. Comparing the polynomials in $t=\ln x$ at each power of $x$ on the both sides, first we have
$$
L\Bigl(1+\frac{d}{dt}\Bigr)P_1(t)=M(0,t,0,\dots,0).
$$
This is an inhomogeneous linear ODE with constant coefficients, with respect to the unknown $P_1$.
Since $L(1)\ne0$ by Lemma 1, zero is not a root of the characteristic polynomial of this equation. Therefore
this ODE has a unique polynomial solution and the degree of this polynomial coincides with that of
the right hand side:
$$
\deg P_1(t)=\deg M(0,t,0,\dots,0)\leqslant C.
$$

Let us denote by $P_k^j(t)$ the polynomial $\bigl(k+\frac{d}{dt}\bigr)^jP_k(t)$, $j=0,1,\ldots,n$ (in particular, $P_k^0=P_k$). Then
$$
\delta^j\psi=\sum_{k=1}^{\infty}P_k^j(\ln x)\,x^k.
$$
Suppose $M$ is a linear combination of monomials of the form $x^{\mu}(\ln x)^{\nu}u^{q_0}(\delta u)^{q_1}\ldots(\delta^nu)^{q_n}$.
Consequently, for each $P_k(t)$ with $k\geqslant2$ we have an inhomogeneous linear ODE with constant coefficients,
\begin{equation}\label{lODE}
L\Bigl(k+\frac{d}{dt}\Bigr)P_k(t)=R_k(t),
\end{equation}
where $R_k(t)$ is a linear combination of polynomials of the form
$$
t^{\nu}\,(P^0_{k_1}\ldots P^0_{k_{q_0}})(P^1_{l_1}\ldots P^1_{l_{q_1}})\ldots(P^n_{m_1}\ldots P^n_{m_{q_n}}),
$$
furthermore
$$
\nu\leqslant C \qquad {\rm and} \qquad \sum_{i=1}^{q_0}k_i+\sum_{i=1}^{q_1}l_i+\ldots+\sum_{i=1}^{q_n}m_i\leqslant k-1.
$$
In view of the inductive assumption,
\begin{eqnarray*}
\deg P^0_{k_1}\ldots P^0_{k_{q_0}} & \leqslant & (k_1+\ldots+k_{q_0})C, \\
\deg P^1_{l_1}\ldots P^1_{l_{q_1}} & \leqslant & (l_1+\ldots+l_{q_1})C, \\
\ldots & \ldots & \ldots \\
\deg P^n_{m_1}\ldots P^n_{m_{q_n}} & \leqslant & (m_1+\ldots+m_{q_n})C,
\end{eqnarray*}
hence
$$
\deg R_k(t)\leqslant C+(k-1)C=k\,C.
$$
Again, since $L(k)\ne0$ by Lemma 1, zero is not a root of the characteristic polynomial of the linear ODE (\ref{lODE}). Therefore
this ODE has a unique polynomial solution $P_k$ whose degree coincides with that of the right hand side:
$$
\deg P_k(t)=\deg R_k(t)\leqslant k\,C.
$$
The lemma is proved.{\hfill $\Box$}

\section{From ODEs to Linear Algebra}

Let us rewrite $\psi$ in the form
$$\psi=\sum\limits_{k=1}^{\infty}P_k(-\epsilon\ln x)\,x^k,$$
where a small constant $\epsilon>0$ will be chosen later (we denote new polynomials by $P_k$ again).
Then the operators $\delta$ and $L(\delta)$ act on a monomial $P_k(-\epsilon\ln x)\,x^k$ as follows:
\begin{eqnarray*}
\delta: & & P_k(-\epsilon\ln x)\,x^k\;\mapsto\;x^k\,\Bigl(k-\epsilon\,\frac{d}{dt}\Bigr)P_k(t)|_{t=-\epsilon\ln x},\\
L(\delta): & & P_k(-\epsilon\ln x)\,x^k\;\mapsto\;x^k\,L\Bigl(k-\epsilon\,\frac{d}{dt}\Bigr)P_k(t)|_{t=-\epsilon\ln x}.
\end{eqnarray*}

This action is naturally represented on the level of vectors and matrices: if ${\bf b}_k$ is a column of the coefficients of the polynomial
$P_k$ and ${\bf c}_k$, ${\bf d}_k$ are columns of the coefficients of the polynomials $\bigl(k-\epsilon\,\frac{d}{dt}\bigr)P_k$,
$L\bigl(k-\epsilon\,\frac{d}{dt}\bigr)P_k$ respectively, then
$$
{\bf c}_k=(k\,I-N_k)\,{\bf b}_k, \qquad {\bf d}_k=L(k\,I-N_k)\,{\bf b}_k,
$$
where  $I$ is an identity matrix, $N_k$ is a nilpotent matrix of the form
$$
N_k=\begin{pmatrix} 0 & \epsilon & 0 & \ldots & 0 & 0 \\ 0 & 0 & 2\epsilon & 0 & \ldots & 0 \\
\hdotsfor{6} \\ \hdotsfor{6} \\ 0 & \ldots & \ldots & 0 & 0 & \nu_k\,\epsilon \\
0 & \ldots & \ldots & \ldots & 0 & 0
\end{pmatrix}, \qquad N_k^{\nu_k+1}=0.
$$

Let us decompose the polynomial $L(\xi)=a_0+\ldots+a_n(\xi+\ell)^n=a_n\prod\limits_{j=1}^n(\xi+\lambda_j)$, ${\rm Re}\,\lambda_j\geqslant 0$. Then
the matrix $L(k\,I-N_k)$ is represented in the form
\begin{equation*}
L(k\,I-N_k)=a_n\prod\limits_{j=1}^n \Bigl((k+\lambda_j)I-N_k\Bigr)=a_n\prod\limits_{j=1}^n(k+\lambda_j)\,\prod\limits_{j=1}^n\Bigl(I-\frac{N_k}{k+\lambda_j}\Bigr)
=L(k)\,\prod\limits_{j=1}^n\Bigl(I-\frac{N_k}{k+\lambda_j}\Bigr).
\end{equation*}
The inverse matrix has the form
$$
L(k\,I-N_k)^{-1}=\frac{1}{L(k)}\,\prod\limits_{j=1}^n \Bigl(I-\frac{N_k}{k+\lambda_j}\Bigr)^{-1},
$$
where
$$
\Bigl(I-\frac{N_k}{k+\lambda_j}\Bigr)^{-1}=I+\frac{N_k}{k+\lambda_j}+\Bigl(\frac{N_k}{k+\lambda_j}\Bigr)^2+\dots+\Bigl(\frac{N_k}{k+\lambda_j}\Bigr)^{\nu_k}.
$$

Further for matrices we will use the matrix 1-norm $\|A\|=\|A\|_1=\max\limits_{j}\sum\limits_{i}|a_{ij}|$ corresponding to the vector 1-norm $\|x\|_1=\sum\limits_{i}|x_i|$.
\medskip

{\bf Lemma 3.} {\it If $\epsilon$ is small enough then there is $0<\varepsilon<1$ such that
\begin{eqnarray*}
\|L(k\,I-N_k)\| & \leqslant & (1+\varepsilon)^n\,|L(k)|, \\
\bigl\|L(k\,I-N_k)^{-1}\bigr\| & \leqslant & \frac1{(1-\varepsilon)^n\,|L(k)|},
\end{eqnarray*}
for any polynomial $L$ of degree $n$ not vanishing in the open right-half plane. In particular,
$\|k\,I-N_k\|\leqslant (1+\varepsilon)k$ and $\bigl\|(k\,I-N_k)^{-1}\bigr\|\leqslant\frac1{(1-\varepsilon)k}$.
}
\medskip

{\bf Proof.} If $\epsilon$ is small enough then $\|N_k\|=\nu_k\,\epsilon\leqslant\varepsilon\,k$, with $\varepsilon$ small enough. Hence, for any integer $k>0$ one has
$$
\Bigl\|I-\frac{N_k}{k+\lambda_j}\Bigr\|\leqslant 1+\varepsilon, \qquad j=1,\dots, n.
$$
Therefore,
$$
\|L(k\,I-N_k)\|=\Bigl\|L(k)\prod\limits_{j=1}^n\Bigl(I-\frac{N_k}{k+\lambda_j}\Bigr)\Bigr\|\leqslant
|L(k)|\,\prod\limits_{j=1}^n\Bigl\| \Bigl(I-\frac{N_k}{k+\lambda_j}\Bigr)\Bigr\| \leqslant |L(k)|\,(1+\varepsilon)^n.
$$
To estimate the norm of the inverse matrix, we note that
$$
\Bigl\|\Bigl(I-\frac{N_k}{k+\lambda_j}\Bigr)^{-1}\Bigr\| \leqslant \|I\|+\Bigl\|\frac{N_k}{k+\lambda_j}\Bigr\|+\dots+
\Bigl\|\frac{N_k}{k+\lambda_j}\Bigr\|^{\nu_k} \leqslant 1+\varepsilon+\dots+\varepsilon^{\nu_k}\leqslant \frac{1}{1-\varepsilon},
$$
whence
$$
\bigl\|L(k\,I-N_k)^{-1}\bigr\|=\Bigl\|\frac{1}{L(k)}\,\prod\limits_{j=1}^n\Bigl(I-\frac{N_k}{k+\lambda_j}\Bigr)^{-1}\Bigr\| \leqslant
\frac{1}{|L(k)|\,(1-\varepsilon)^n}.
$$
The lemma is proved. {\hfill $\Box$}

\section{A majorant equation}

Let us rewrite the reduced equation (\ref{rODE}) in the form
\begin{equation}\label{ODE2}
L(\delta)u=x\,M(x,-\epsilon\ln x,u,\delta u,\dots,\delta^n u),
\end{equation}
denoting the new polynomial on the right hand side by the same letter $M$.
Consider another equation which, as we will see later, is a majorant one for (\ref{ODE2}) in some sense:
\begin{equation}\label{ODE3}
\sigma\,\delta^nu=x\,\widetilde{M}(x,-\epsilon\ln x,\delta^nu),
\end{equation}
where
$$
\frac1{\sigma}:=\sup\limits_{k\geqslant1}\,\bigl\|L(k\,I-N_k)^{-1}\bigr\|\cdot\bigl\|(k\,I-N_k)^n\bigr\|<+\infty
$$
by Lemma 3, and a polynomial $\widetilde{M}$ is constructed as follows. Suppose $M$ is a sum of monomials of the form
\begin{equation}\label{Msummand}
\alpha\,x^{\mu}(-\epsilon\ln x)^{\nu}u^{q_0}(\delta u)^{q_1}\ldots(\delta^nu)^{q_n}, \quad \alpha\in\mathbb C.
\end{equation}
To consrtuct the polynomial $\widetilde{M}$, one changes each such a summand to
\begin{equation}\label{tMsummand}
|\alpha|\,x^{\mu}(-\epsilon\ln x)^{\nu}(c\,\delta^nu)^{q_0}(c\,\delta^nu)^{q_1}\ldots(c\,\delta^nu)^{q_n}, \qquad
c=\Bigl(\frac{1+\varepsilon}{1-\varepsilon}\Bigr)^n.
\end{equation}
\medskip

{\bf Lemma 4.} {\it The equation $(\ref{ODE3})$ possesses a uniquely determined formal Dulac series solution
$$
\tilde{\psi}=\sum\limits_{k=1}^{\infty}\widetilde P_k(-\epsilon\ln x)\,x^k,
$$
where $\widetilde P_k\in\mathbb{R}_+[t]$ are polynomials with non-negative real coefficients of degree
$\deg\widetilde P_k\leqslant k\,C$.
}
\medskip

{\bf Proof.} A proof is analogous to that of Lemma 2: each $\widetilde P_k$ is obtained as a solution of an inhomogeneous linear ODE
with constant coefficients. Starting with
$$
\sigma\Bigl(1-\epsilon\,\frac{d}{dt}\Bigr)^n\,\widetilde P_1(t)=\widetilde{M}(0,t,0)\in\mathbb{R}_+[t],
$$
then one finds the other $\widetilde P_k$, $k\geqslant2$, as unique polynomial solutions of the corresponding ODE
$$
\sigma\Bigl(k-\epsilon\,\frac{d}{dt}\Bigr)^n\,\widetilde P_k(t)=\widetilde Q_k(t)\in\mathbb{R}_+[t], \qquad k\geqslant2,
$$
where the polynomial $\widetilde Q_k$ is expressed {\it via} the polynomials $\widetilde Q_1=\widetilde{M}(0,t,0),
\widetilde Q_2,\ldots,\widetilde Q_{k-1}$ and its degree does not exceed $k\,C$ (see the detailed expressions in the proof
of the next Lemma 5). Thus, the non-negativity of the coefficients of each polynomial $\widetilde P_k$
follows from the non-negativity of the coefficients of the corresponding $\widetilde Q_k$ and from the non-negativity
of the elements of the matrix $(k\,I-N_k)^{-1}$. {\hfill $\Box$}
\medskip

For an arbitrary polynomial $P\in{\mathbb C}[t]$, let us define its norm $\|P\|$ as the 1-norm of the column of its coefficients.
Besides the standard norm properties, this satisfies the following ones which are not difficult to check:

1) for any $P, Q\in{\mathbb C}[t]$ one has $\|PQ\|\leqslant\|P\|\cdot\|Q\|$;

2) if $P, Q\in{\mathbb R}_+[t]$ then $\|P+Q\|=\|P\|+\|Q\|$ and $\|PQ\|=\|P\|\cdot\|Q\|$.
\medskip

Now we will prove that the constructed equation (\ref{ODE3}) is a majorant one for the initial equation (\ref{ODE2}) in the following sense.
\medskip

{\bf Lemma 5.} {\it The formal Dulac series solution $\tilde\psi$ of the equation $(\ref{ODE3})$ is majorant for the formal Dulac series solution
$\psi$ of the equation $(\ref{ODE2}):$ $\|P_k\|\leqslant\| \widetilde P_k\|$ for all $k$.}
\medskip

{\bf Proof.} We have already understood that the polynomials $P_k$ and $\widetilde P_k$ are solutions of the corresponding
inhomogeneous linear ODEs with constant coefficients:
\begin{eqnarray}\label{ODEforP}
L\Bigl(k-\epsilon\,\frac{d}{dt}\Bigr)P_k(t)&=&Q_k(t), \nonumber \\
\sigma\Bigl(k-\epsilon\,\frac{d}{dt}\Bigr)^n\widetilde P_k(t)&=&\widetilde Q_k(t),
\end{eqnarray}
where $Q_1(t)=M(0,t,0,\ldots,0)$ and $\widetilde Q_1(t)=\widetilde{M}(0,t,0)$, whence $\|Q_1\|=\|\widetilde Q_1\|$. Therefore
\begin{eqnarray*}
\|P_1\|\leqslant\bigl\|L(I-N_1)^{-1}\bigr\|\cdot\|Q_1\|=\bigl\|L(I-N_1)^{-1}\bigr\|\cdot\|\widetilde Q_1\|\leqslant \\
\leqslant\sigma\bigl\|L(I-N_1)^{-1}\bigr\|\cdot\bigl\|(I-N_1)^n\bigr\|\cdot\|\widetilde P_1\|\leqslant\|\widetilde P_1\|.
\end{eqnarray*}

To obtain analogous estimates for all the other $k\geqslant2$, let us study expressions for the corresponding $Q_k$ and $\widetilde Q_k$ in more details. We will use the notation $P_k^j(t)$ for the polynomial
$\bigl(k-\epsilon\,\frac{d}{dt}\bigr)^jP_k(t)$, $j=0,1,\ldots,n$ (in particular, $P_k^0=P_k$). Then
\begin{eqnarray*}
\delta^j\psi=\sum_{k=1}^{\infty}P_k^j(-\epsilon\ln x)\,x^k,\\
\delta^n\tilde\psi=\sum_{k=1}^{\infty}\frac1{\sigma}\,\widetilde Q_k(-\epsilon\ln x)\,x^k.
\end{eqnarray*}
Looking at (\ref{Msummand}) and (\ref{tMsummand}) one concludes that $Q_k(t)$ is a sum of polynomials of the form
\begin{equation}\label{Qsummand}
\alpha\,t^{\nu}\,(P^0_{k_1}\ldots P^0_{k_{q_0}})(P^1_{l_1}\ldots P^1_{l_{q_1}})\ldots(P^n_{m_1}\ldots P^n_{m_{q_n}}),
\end{equation}
where $\sum_{i=1}^{q_0}k_i+\sum_{i=1}^{q_1}l_i+\ldots+\sum_{i=1}^{q_n}m_i\leqslant k-1$, and $\widetilde Q_k(t)$ is
a sum of the corresponding polynomials
\begin{equation}\label{tQsummand}
|\alpha|\,t^{\nu}\,\Bigl(\frac c{\sigma}\,\widetilde Q_{k_1}\ldots\frac c{\sigma}\,\widetilde Q_{k_{q_0}}\Bigr)
\Bigl(\frac c{\sigma}\,\widetilde Q_{l_1}\ldots\frac c{\sigma}\,\widetilde Q_{l_{q_1}}\Bigr)\ldots
\Bigl(\frac c{\sigma}\,\widetilde Q_{m_1}\ldots\frac c{\sigma}\,\widetilde Q_{m_{q_n}}\Bigr).
\end{equation}
The norm of (\ref{Qsummand}) does not exceed
$$
|\alpha|\cdot\|P^0_{k_1}\|\ldots\|P^0_{k_{q_0}}\|\cdot\|P^1_{l_1}\|\ldots\|P^1_{l_{q_1}}\|\cdot\ldots\cdot\|P^n_{m_1}\|\ldots\|P^n_{m_{q_n}}\|.
$$
Using the inductive assumption and the relations (\ref{ODEforP}) we can estimate each factor $\|P_s^j\|$, $s<k$, in the product above:
\begin{eqnarray*}
\|P_s^j\| & \leqslant & \bigl\|(s\,I-N_s)^j\bigr\|\cdot\|P_s\|\leqslant\bigl\|(s\,I-N_s)^j\bigr\|\cdot\|\widetilde P_s\| \leqslant \\
 & \leqslant & \bigl\|(s\,I-N_s)^j\bigr\|\cdot\bigl\|(s\,I-N_s)^{-n}\bigr\|\cdot\|\widetilde Q_s/\sigma\| \leqslant
 \|s\,I-N_s\|^j\cdot\bigl\|(s\,I-N_s)^{-1}\bigr\|^n\cdot\|\widetilde Q_s/\sigma\| \leqslant \\
 & \leqslant & \frac{(1+\varepsilon)^js^j}{(1-\varepsilon)^ns^n}\,\|\widetilde Q_s/\sigma\|\leqslant\frac c{\sigma}\,
\|\widetilde Q_s\|.
\end{eqnarray*}
Hence, the norm of (\ref{Qsummand}) does not exceed the norm of (\ref{tQsummand}) (recall that the norm of a product equals the product of norms for polynomials with non-negative real coefficients) and, therefore,
$\|Q_k\|\leqslant\|\widetilde Q_k\|$ (again, the norm of a sum equals the sum of norms
for polynomials with non-negative real coefficients).

Finally, we conclude
\begin{eqnarray*}
\|P_k\|\leqslant\bigl\|L(k\,I-N_k)^{-1}\bigr\|\cdot\|Q_k\|\leqslant\bigl\|L(k\,I-N_k)^{-1}\bigr\|\cdot\|\widetilde Q_k\|
\leqslant \\
\leqslant\sigma\bigl\|L(k\,I-N_k)^{-1}\bigr\|\cdot\bigl\|(k\,I-N_k)^n\bigr\|\cdot\|\widetilde P_k\|\leqslant\|\widetilde P_k\|.
\end{eqnarray*}
The lemma is proved. {\hfill $\Box$}

\section{Proof of the theorem on the convergence}

Since the function $\ln x$ is transcendental, the majorant equation \eqref{ODE3}
can be regarded as an algebraic one,
\begin{equation}\label{algebraic}
\sigma\,U=x\,\widetilde{M}(x,t,U)
\end{equation}
(with two independent variables $x, t$ and the unknown $U=\delta^nu$)  having a formal solution
$$
\widehat U=\sum\limits_{k=1}^{\infty}\frac1{\sigma}\,\widetilde Q_k(t)\,x^k, \qquad
\widetilde Q_k\in\mathbb{R}_+[t].
$$
After opening the brackets in the series $\widehat U$ we get the bivariate power series
$\widehat U_{\rm pow}=\sum\limits_{k=1}^{\infty}\sum\limits_{l=0}^{\tilde\nu_k}c_{kl}\,t^lx^k$, $c_{kl}\in\mathbb{R}_+$,
which also formally satisfies the equation (\ref{algebraic}) (to open the brackets in $\widehat U$ and then to substitute
$\widehat U_{\rm pow}$ into the both sides of (\ref{algebraic}) is the same as to substitute $\widehat U$ into
(\ref{algebraic}) at first and then to open the brackets on its both sides). By the implicit function theorem, this series converges absolutely for small $t$ and $x$ but this is not what we finally need, as the variable $t$ responses for $\ln x$ in the Dulac series, and the latter is unbounded for small $x$. Therefore we take an integer $r>C$ and consider an open sector $S$ with the vertex at the origin and of the opening less than $2\pi$, such that
$$
\left|-\epsilon\ln x\right|<|x|^{-1/r} \qquad \forall x\in S,
$$
which implies
\begin{equation}\label{ineq_on_ln2}
\bigl|\widetilde Q_k(-\epsilon\ln x)\bigr|<\widetilde Q_k\bigl(|x|^{-1/r}\bigr) \qquad \forall x\in S.
\end{equation}

Now let us consider an equation
\begin{equation} \label{maj_eq_modif}
\sigma\,U=x\,\widetilde{M}(x,x^{-1/r}, U)
\end{equation}
obtained from (\ref{algebraic}) by putting $t=x^{-1/r}$. This has a formal Puiseux series solution
$$
\phi=\sum\limits_{k=1}^{\infty}\sum\limits_{l=0}^{\tilde\nu_k}c_{kl}\,x^{k-l/r}
$$
obtained, respectively, from the bivariate power series $\widehat U_{\rm pow}$. Indeed, to put $t=x^{-1/r}$
in $\widehat U_{\rm pow}$ and then to substitute the obtained Puiseux series into the both sides of (\ref{maj_eq_modif})
is the same as to substitute $\widehat U_{\rm pow}$ into the equation (\ref{algebraic}) at first and then to put $t=x^{-1/r}$
on its both sides. (Note that $k-l/r>0$ for all $k$, $l$, since $\tilde\nu_k\leqslant k\,C<kr$, and by the same reason
$k_1-l_1/r\ne k_2-l_2/r$ if $(k_1, l_1)\ne(k_2, l_2)$.) The Puiseux series $\phi$ converges in $S$ absolutely for $x$ small enough (to see this, it is sufficient to make a change $x=z^r$ in (\ref{maj_eq_modif}) and apply the implicit function theorem).

The series
$$
\phi^{\circ}(|x|)=\sum\limits_{k=1}^{\infty}\frac1{\sigma}\,\widetilde Q_k(|x|^{-1/r})\,|x|^k
$$
is another representation of the Puiseux series $\phi(|x|)=\sum\limits_{k=1}^{\infty}\sum\limits_{l=0}^{\tilde\nu_k}
c_{kl}\,|x|^{k-l/r}$, hence it also converges in $S$ for $x$ small enough (say, for $|x|<\rho$) by the corresponding well known property of
convergent positive series (see \cite[Ch. VIII]{Gu}). From the inequality \eqref{ineq_on_ln2} it follows that the series
$$
\delta^n\tilde\psi=\sum\limits_{k=1}^{\infty}\frac1{\sigma}\,\widetilde Q_k(-\epsilon\ln x)\,x^k
$$
converges in $S$ absolutely for $|x|<\rho$, therefore $\tilde\psi=\sum_{k=1}^{\infty}\widetilde P_k(-\epsilon\ln x)\,x^k$
does. Since $-\epsilon\ln |x|\leqslant |\epsilon\ln x|$, we then have the convergence of the series
$$
\sum\limits_{k=1}^{\infty}\widetilde P_k(-\epsilon\ln |x|)\,|x|^k, \qquad |x|<\rho.
$$

Now we finish the proof of Theorem 1 by proving the convergence of the series $\psi=\sum_{k=1}^{\infty}
P_k(-\epsilon\ln x)\,x^k$. We have
$$
\bigl|P_k(-\epsilon\ln x)\,x^k\bigr|\leqslant\|P_k\|\cdot|\epsilon\ln x|^{\nu_k}\,|x|^k\leqslant\|P_k\|\cdot|x|^{(1-C/r)k}
\leqslant\|\widetilde P_k\|\cdot|x|^{(1-C/r)k}.
$$
Thus, for all $x$ such that $w=|x|^{1-C/r}<\rho$ one has
$$
\bigl|P_k(-\epsilon\ln x)\,x^k\bigr|\leqslant\|\widetilde P_k\|\,w^k\leqslant\widetilde P_k(-\epsilon\ln w)\,w^k,
$$
whence the convergence of $\psi$ follows. The theorem is proved.

\medskip
{\bf Remark.} Using the results of the present work and article \cite{GG} one can prove a statement similar to Theorem 1 but for formal Dulac series of a more general form (a proof is more technical though). {\it Let the series
$$
\varphi=\sum\limits_{k=0}^{\infty}p_k(\ln x)x^{\lambda_k},\qquad \lambda_k\in\mathbb{C},
\qquad 0\leqslant{\rm Re\,}{\lambda_0}\leqslant {\rm Re\,}{\lambda_1}\leqslant\dots\rightarrow +\infty,$$
formally satisfy the equation under consideration:
$$F(x,\Phi):=F(x,\varphi,\delta\varphi,\dots,\delta^n \varphi)=0,$$
and let for each $j=0,\dots,n$ one have
$$
\frac{\partial F}{\partial y_j}(x,\Phi)=a_j x^{\alpha}+b_j(\ln x) x^{\alpha_j}+\dots,
\qquad {\rm Re\,}\alpha<{\rm Re\,}\alpha_j,
$$
where $a_j\in\mathbb{C}$, $\alpha\in\mathbb{C}$ is the same for all $j$, and $b_j\in\mathbb{C}[t]$.
If $a_n\neq 0$ then for any open sector $S$ of sufficiently small radius, with the vertex at the origin and of the opening less than $2\pi$, the series $\varphi$ converges uniformly in $S$.}

\section{Examples}

In this section we give two examples, of the Abel equation and of the Painlev\'e VI equation, and apply Theorem 1 for
proving the convergence of their formal Dulac series solutions.
\medskip

{\bf Example 1.} The Abel equation of the second kind 
$$
w\,\frac{dw}{dx}=-w-\frac{1}{x^3}
$$
has a one-parameter family of formal Dulac series solutions
$$\hat w=\frac{1}{x}\Bigl(1+(C-\ln x)x^2+\sum\limits_{k=2}^{\infty}P_k(\ln x)x^{2k}\Bigr),\qquad C\in\mathbb{C}.
$$
Making the power transformation $w=y/x$ in the equation under consideration and rewriting the result by means of the operator $\delta$, we obtain the equation
$$
f(x,y,\delta y):=y\,\delta y-y^2+x^2y+1=0
$$
with a family of formal Dulac series solutions
$$
\varphi=1+(C-\ln x)x^2+\sum\limits_{k=2}^{\infty}P_k(\ln x)x^{2k}.
$$
Let us check their convergence by applying Theorem 1. As $f(x,y_0,y_1)=y_0y_1-y_0^2+x^2y_0+1$, one has
$$
\frac{\partial f}{\partial y_0}=y_1-2y_0+x^2,\qquad \frac{\partial f}{\partial y_1}=y_0.
$$
Substituting
\begin{eqnarray*}
y_0&=&\varphi=1+(C-\ln{x})x^2+\dots, \\
y_1&=&\delta\varphi=(2C-1-2\ln{x})x^2+\dots,
\end{eqnarray*}
we thus obtain
$$
\frac{\partial f}{\partial y_0}=-2+\dots, \qquad \frac{\partial f}{\partial y_1}=1+(C-\ln{x})x^2+\dots\,.
$$
Hence, the condition of Theorem 1 is fulfilled and the series $\varphi$ converges in any open sector
$S\subset \mathbb{C}$ with the vertex at the origin, of sufficiently small radius and of the opening less than $2\pi$.
\medskip

{\bf Example 2.} The Painlev\'e VI equation 
\begin{equation}\label{painleve6}
y''=\frac{(y')^2}{2}\left(\frac{1}{y}+\frac{1}{y-1}+\frac{1}{y-x}\right)-y'\left(\frac{1}{x}
+\frac{1}{x-1}+\frac{1}{y-x}\right)+
\end{equation}
$$+\frac{y(y-1)(y-x)}{x^2(x-1)^2}\left[a+b\frac{x}{y^2}+c\frac{x-1}{(y-1)^2}+
d\frac{x(x-1)}{(y-x)^2}\right],\qquad a,\;b,\;c,\;d\in\mathbb{C},$$ 
with $a, c\neq 0\;$, $a\neq c$, and $\;\sqrt{2a}\pm\sqrt{2c}\in\mathbb{N}\;$ 
has two one-parameter families of formal Dulac series solutions
$$
\varphi=1\pm\sqrt{c/a}+\sum\limits_{k=1}^\infty P_k(\ln{x})x^k,
$$
where a free parameter is contained in a polynomial $P_k$ with $k=\sqrt{2a}\pm\sqrt{2c}$, see \cite{BG}.
One rewrites equation \eqref{painleve6} in the form
\begin{equation}\label{painleve_rew}
\begin{array}{ccc}f(x,\,y_0,\,y_1,\,y_2):=
2y_2(x-1)^2y_0(y_0-1)(y_0-x)-y_1^2(x-1)^2[(y_0-1)(y_0-x)+\\
+y_0(y_0-x)+y_0(y_0-1)\;]+2y_1 x(x-1)y_0(y_0-1)^2-\\
-[\;2ay_0^2(y_0-1)^2(y_0-x)^2+2bx(y_0-1)^2(y_0-x)^2+\\
+2c(x-1)y_0^2(y_0-x)^2+2dx(x-1)y_0^2(y_0-1)^2\,]=0.\end{array}
\end{equation}
The partial derivative of $f$ with respect to $y_2$ is
$$
\frac{\partial f}{\partial y_2}=2(x-1)^2y_0(y_0-1)(y_0-x).
$$
Substituting $y_0=\varphi$ in this expression we obtain
$$
\frac{\partial f}{\partial y_2}=\pm 2\left(1\pm\sqrt{c/a}\right)^2 \sqrt{c/a}+b_2(\ln{x})x+\dots,
$$
for some polynomial $b_2$. As $c\neq0$ and $a\neq c$, this Dulac series begins with a non-zero constant, which is sufficient
for the condition of Theorem 1 being fulfilled. Thus, the series $\varphi$ converges in any open sector $S\subset \mathbb{C}$ 
with the vertex at the origin, of sufficiently small radius and of the opening less than $2\pi$.

As we mentioned in Introduction, in such a way one can prove the convergence of all formal Dulac series solutions of the third, fifth and sixth 
Painlev\'e equations (which has been also done by S.\,Shimomura for the Painlev\'e V and VI in \cite{Sh2} and \cite{Sh1}, respectively).

\end{document}